\documentclass[12pt]{amsart}
\usepackage{amsmath, amsfonts, amssymb, amsthm}

\font\teneufm=eufm10 \font\seveneufm=eufm7 \font\fiveeufm=eufm5
\newfam\eufmfam
\textfont\eufmfam=\teneufm \scriptfont\eufmfam=\seveneufm
\scriptscriptfont\eufmfam=\fiveeufm

\newfam\msbfam
\font\tenmsb=msbm10 scaled \magstep1   \textfont\msbfam=\tenmsb
\font\sevenmsb=msbm7  scaled \magstep1
\scriptfont\msbfam=\sevenmsb \font\fivemsb=msbm5  scaled \magstep1
\scriptscriptfont\msbfam=\fivemsb
\def\Bbb{\fam\msbfam \tenmsb}

\def\BB{{\Bbb B}}

\def\CC{{\Bbb C}}

\newfam\bobfam
\font\tenbob=msbm10    \textfont\bobfam=\tenbob
\font\sevenbob=msbm7     \scriptfont\bobfam=\sevenbob
\font\fivebob=msbm5     \scriptscriptfont\bobfam=\fivebob

 \def\HollowBox #1#2{{\dimen0=#1 \advance\dimen0 by -#2
       \dimen1=#1 \advance\dimen1 by #2
        \vrule height #1 depth #2 width #2
        \vrule height 0pt depth #2 width #1
        \llap{\vrule height #1 depth -\dimen0 width \dimen1}%
       \hskip -#2
       \vrule height #1 depth #2 width #2}}

\def\endpf{\hfill $\BoxOpTwo$}

\def\endpf{\hfill $\Box$}

\def\Re{{\rm Re}\,}

\def\bo{\partial \Omega}

\def\bo{\partial \Omega}

\newtheorem{theorem}{Theorem}[section]
\newtheorem{lemma}{Lemma}[section]
\newtheorem{proposition}{Proposition}[section]

\newtheorem{remark}{Remark}[section]
\newtheorem{definition}{Definition}[section]
\newtheorem{question}{Question}[section]

\begin{document}

\title{A Kobayashi metric version of Bun Wong's Theorem}

\author{Kang-Tae Kim and Steven G. Krantz}

\thanks{The research of the first named author is supported
in part by the KOSEF Grant R01-2005-000-10771-0 of The Republic of
Korea. He also would like to thank the Centre de Mathematiques et
Informatiques, Universit\'e de Provence, Aix-Marseille 1 (France)
for its hospitality at the final stage of this work.}

\thanks{Both authors thank the American Institute of
Mathematics for its hospitality during a portion of this work.}

\thanks{The second named author has grants from the National
Science Foundation and the Dean of the Graduate School at
Washington University.}

\maketitle

\parindent=.22in

\section{Introduction}

\subsection{Basic Terminology and Statement of Main Theorem}

For a complex manifold $M$, denote by $k_M$ its
Kobayashi-Royden infinitesimal metric and by $d_M$ its
Kobayashi distance, and (see [KRA1], [KOB2]).

\begin{definition} \rm
A map $f:M \to N$ from a complex manifold $M$ into another complex
manifold $N$ is said to be a {\it Kobayashi isometry} if $f$ is a
homeomorphism satisfying the condition that $d_N (f(x), f(y)) =
d_M (x,y)$ for every $x,y \in M$.
\end{definition}

The set $G_M$ of Kobayashi isometries of $M$ (onto $M$ itself)
endowed with the compact-open topology is a topological group with
respect to the binary law of composition of mappings. We call this group
$G_M$, the {\it Kobayashi isometry group} of the complex manifold
$M$.
\bigskip

Denote by $\BB^n$ the open unit ball in $\CC^n$. Notice that the
Kobayashi distance of $\BB^n$ in fact coincides with the
Poincar\'e-Bergman distance of $\BB^n$.  The primary aim of this
article is to establish the following theorem, which is a Kobayashi
metric version of Bun Wong's classical theorem \cite{WON}.

\begin{theorem} \label{main}
Let $\Omega$ be a bounded domain in $\CC^n$ with a
$\mathcal{C}^{2,\epsilon}$ smooth ($\epsilon>0$), strongly
pseudoconvex boundary. If its Kobayashi isometry group $G_\Omega$
is non-compact, then $\Omega$ is biholomorphic to the open unit
ball $\BB^n$.
\end{theorem}

\begin{remark} \rm
The classical results (see [WON] and [ROS]) assume noncompact
(biholomorphic) automorphism group of the domain $\Omega$.
But it has been understood for many years that this
theorem of Bun Wong and Rosay is really a ``flattening'' result
of geometry (terminology of Gromov is being used here).  Our new
theorem puts this relationship into perspective.
\endpf
\end{remark}

\subsection{The Kobayashi-distance version of Wong's theorem}
Experts who are familiar with Wong's theorem \cite{WON} would
expect the following:

\begin{theorem}[Seshadri-Verma] \label{main-1}
Let $\Omega$ be a bounded domain in $\CC^n$ with a $\mathcal{C}^2$
smooth, strongly pseudoconvex boundary. If its Kobayashi isometry
group $G_\Omega$ is non-compact, then there exists a Kobayashi
isometry $f:\Omega \to \BB^n$.
\end{theorem}

Certainly it was Seshadri and Verma ([SEV1, SEV2] as well as
[VER]) who first conceived the idea of a metric version of
Wong's theorem. We present here a different proof of Theorem
1.2 which closes some gaps in the extant argument and answers
some subtle questions. We believe that the arguments we
present in this article can be of use for many other purposes.
The arguments involved here are subtle, because the mappings
under consideration are {\sl a priori} only continuous. Unlike
the holomorphic case the restrictions of Kobayashi isometries
to sub-domains are not isometries with respect to the
Kobayashi metric of the sub-domain. Another slippery point is
that the full power of Montel's theorem and Cartan's
uniqueness theorem is not available for equi-continuous maps.
So it is necessary to give precise estimates and arguments
that clarify all these subtle points necessary for the proof.
It is true that our proof of this theorem follows
the same general line of reasoning as \cite{SEV1, SEV2}, which
in turn is along the line of scaling method introduced by S.
Pinchuk around 1980 (\cite{PIN}). \medskip

The second half, that is indeed the main part of this
article, presents the following:

\begin{theorem} \label{main-main}
Let $\Omega$ be a bounded domain in $\CC^n$ with a
$\mathcal{C}^{2,\epsilon}$ smooth, strongly pseudoconvex boundary.
If there exists a Kobayashi isometry $f:\Omega \to \BB^n$ from
$\Omega$ onto the open unit ball $\BB^n$ of $\CC^n$, then $f$ is
either holomorphic or conjugate holomorphic.
\end{theorem}

Notice that Theorems \ref{main-1} and \ref{main-main} imply
Theorem \ref{main}.  And that is the main point of this paper.
The work [SEV2], which addresses similar questions, assumes that
the mapping is ${\mathcal C}^1$ to the boundary.  We are able
to eliminate this somewhat indelicate hypothesis.
\medskip

Seshadri and Verma in the above cited work have proved the same
conclusion in the case when $\Omega$ is strongly convex. Notice
that our theorem here assumes only {\it strong pseudoconvexity}.
On the other hand, the experts in this line of research would feel
that the optimal regularity of the boundary should be
$\mathcal{C}^2$, instead of $\mathcal{C}^{2,\epsilon}$ with some
$\epsilon>0$. At the time of this writing, we do not know how to
achieve the optimum, because of technical reasons connected with
work of Lempert \cite{LEM1, LEM2}. We would like to mention it as
a question for future study.

\section{Some Fundamentals}

\subsection{Terminology and Notation}
Let $\Omega$ be a Kobayashi hyperbolic domain in $\CC^n$. For a
point $q \in \Omega$, let us write
$$
G_\Omega (q) = \{ \varphi (q) \mid \varphi \in G_\Omega \}.
$$
This is usually called the (point) {\it orbit} of $q$ under the action of
the Kobayashi isometry group $G_\Omega$ on the domain $\Omega$.
\bigskip

Call a boundary point $p \in \partial\Omega$ a {\it boundary orbit
accumulation point} if $p$ belongs to the closure
$\overline{G_\Omega (q)}$ of the orbit $G_\Omega (q)$ of a certain
interior point $q \in \Omega$ under the action of the Kobayashi
isometry group $G_\Omega$. In other words, $p$ is a boundary orbit
accumulation point if and only if there exists an
interior point $q \in \Omega$ and a sequence of Kobayashi
isometries $\varphi_j \in G_\Omega$ such that
$\displaystyle{\lim_{j \to \infty} \varphi_j (q) = p}$.
\bigskip

Let us adopt the notation
$$
B_d (q; r) := \{ y \mid d(y,q)<r \}
$$
for any distance $d$ in general. Then one observes the following:

\begin{proposition} \label{prop-1}
Let $\Omega$ be a bounded, complete Kobayashi hyperbolic domain in
$\CC^n$. Then its Kobayashi isometry group $G_\Omega$ is
non-compact if and only if $\Omega$ admits a boundary orbit
accumulation point.
\end{proposition}

\it Proof. \rm Notice that the {\em sufficiency} is obvious. We
establish the {\em necessity} only.
\medskip

Expecting a contradiction, assume to the contrary that there are no
boundary orbit accumulation points. Then, for every point $q$ of
the domain $\Omega$, the orbit of $q$ under the group action is
relatively compact. Now let $\{\varphi_j\}$ be an arbitrarily
chosen sequence of Kobayashi isometries; then it is obviously an
equi-continuous family with respect to the Kobayashi distance. By
Barth's theorem (\cite{BAR}), this implies that $\varphi_j$ forms
an equi-continuous family on compact subsets with respect to the
Euclidean distance. Thus one may use the Arzela-Ascoli theorem to
extract a sequence $\{\varphi_{j_k}\}$ that converges uniformly on
compact subsets of $\Omega$ to a limit mapping $\widehat \varphi$.
Thus, replacing $\varphi_j$ by a subsequence, one may assume
without loss of generality that $\varphi_j$ converges uniformly on
compacta to a continuous map, say $\widehat\varphi$.
\medskip

Since the point orbit is always compact, $\widehat\varphi (a) = b$
for some $a, b \in \Omega$. Notice that, exploiting the
completeness of $d_\Omega$, one can deduce that
$$
\widehat\varphi (\Omega) = \widehat\varphi \Big(
\bigcup_{\nu=1}^\infty B_{d_\Omega} (a;\nu) \Big) =
\bigcup_{\nu=1}^\infty \widehat\varphi (B_{d_\Omega} (a;\nu))
\subset \bigcup_{\nu=1}^\infty B_{d_\Omega} (b;\nu) = \Omega.
$$
It is obvious that one can apply the same argument to the sequence
$\varphi_j^{-1}$ (replacing it by a subsequence that converges
uniformly on compacta, if necessary). Thus $\widehat\varphi:\Omega
\to \Omega$ is a homeomorphism. It is obvious that
$\widehat\varphi$ preserves the Kobayashi distance $d_\Omega$.
Altogether, it follows that $\widehat\varphi \in G_\Omega$. This
establishes that $G_\Omega$ is compact. (This argument shows the
sequential compactness of $G_\Omega$, to be precise.  But then
$G_\Omega$ equipped with the topology of uniform convergence on
compacta is metrizable.) This contradiction yields the desired
conclusion.
\endpf
\smallskip \\

Since every bounded strongly pseudoconvex domain is complete with
respect to the Kobayashi distance (see [GRA]), Theorem \ref{main-1} now
follows by the following more general statement:

\begin{theorem} \label{main-2}
Let $\Omega$ be a domain in $\CC^n$. If $\Omega$ admits a boundary orbit
accumulation point at which the boundary of $\Omega$ is
$\mathcal{C}^2$ smooth, strongly pseudoconvex, then $\Omega$ is
Kobayashi isometric to the open unit ball in $\CC^n$.
\end{theorem}

Notice that for this theorem one does not need to assume that the
domain has to be {\it a priori} bounded or Kobayashi hyperbolic.
Instead, the complete Kobayashi hyperbolicity of the domain will
be obtained along the way in the proof, from the given hypothesis
only.  The proof of this result will be developed in Section 3.

\subsection{Some Comparison Estimates}
We now give some comparison inequalities for the Kobayashi
distances and the Kobayashi-Royden infinitesimal metrics for a
sub-domain against its ambient domain. This will play a crucial
role in the proof of Theorem \ref{main-2}.

\begin{lemma}{\rm (Kim-Ma)} \label{lemma-1}
Let $\Omega$ be a Kobayashi hyperbolic domain in $\CC^n$ with a
subdomain $\Omega' \subset \Omega$. Let $q, x \in \Omega'$, let
$d_\Omega (q,x)=a$, and let $b>a$. If $\Omega'$ satisfies the
condition $B_{d_\Omega} (q; b) \subset \Omega'$, then the
following two inequalities hold:
$$
d_{\Omega'} (q,x) \le \frac{1}{\tanh(b-a)} \ d_\Omega (q,x),
$$
$$
k_{\Omega'} (x, v) \le \frac{1}{\tanh(b-a)}~ k_\Omega (x, v),  ~~v
\in \CC^n.
$$
[Recall here that $k$ is the infinitesimal Kobayashi/Royden metric and
$d$ is the integrated Kobayashi/Royden distance.]
\end{lemma}

\it Proof. \rm For the sake of the reader's convenience, we include
here the proof, lifting it from [KIMA]. Let $s=\tanh(b-a)$ and let
$\epsilon > 0$. Denote by $\Delta$ the open unit disc in $\CC$ and
by $\Delta(a;r)$ the open disc of Euclidean radius $r$ centered at $a$ in
$\CC$. Then, by definition of $k_\Omega (x, v)$, there exists a
holomorphic map $h:\Delta \to \Omega$ such that $h(0) = x$ and
$h'(0) = v/(k_\Omega (x,v) + \epsilon)$. If $\zeta \in \Delta
(0;s)$, then
\begin{eqnarray*}
d_\Omega (q, h(\zeta)) %
& \le & d_\Omega (q, x) + d_\Omega (x,h(\zeta)) \\
& = & a + d_\Omega (h(0), h(\zeta)) \\
& \le & a + d_\Delta (0,\zeta) \\
& < & a + (b-a) \\
& = & b.
\end{eqnarray*}
This shows that $h(\Delta (0;s)) \subset \Omega'$.  Now define
$g:\Delta \to \Omega'$ by $g(z) := h(sz)$. Then one has $g(0)=x$
and $g'(0) = sh'(0) = sv/(k_\Omega(x,v) + \epsilon)$. This implies
that $k_{\Omega'} (x,v) \le (k_\Omega (x,v) + \epsilon)/s$. Since
$\epsilon$ is an arbitrarily chosen positive number, it follows
that
$$
k_{\Omega'} (x,v) \le \frac1{\tanh(b-a)}~ k_\Omega (x,v), \quad
\forall v \in \CC^n.
$$

Now let $\delta$ be chosen such that $0 < \delta < b-a$.  There is
a $\mathcal{C}^1$ curve $\gamma : [0,1] \to \Omega$ such that
$\gamma (0) = q$, $\gamma (1) = x$, and $\displaystyle{\int_0^1
k_\Omega (\gamma(t), \gamma'(t))~dt < a+\delta}$. This implies
that $d_\Omega (q, \gamma(t)) < a + \delta < b$ for any $t \in
[0,1]$. Hence $\gamma (t) \in \Omega'$ for every $t \in [0,1]$.
Notice that the inequality
$$
k_{\Omega'} (\gamma(t), \gamma'(t)) \le k_\Omega (\gamma(t),
\gamma'(t))/\tanh (b-a-\delta)
$$
holds for every $t \in [0,1]$ by the preceding arguments. But then
this implies that
\begin{multline*}
d_{\Omega'} (x,q) \le \int_0^1 k_{\Omega'} (\gamma(t), \gamma'(t)) \ dt \\
\le \frac1{\tanh(b-a-\delta)} \int_0^1 k_\Omega (\gamma(t),
\gamma'(t) \ dt.
\end{multline*}
Consequently, one deduces that $d_{\Omega'} (x,q) <
(a+\delta)/\tanh(b-a-\delta)$. Now, letting $\delta$ tend to $0$,
one obtains the desired conclusion.
\endpf
\bigskip \\
\bigskip

\section{Scaling With a Sequence of Kobayashi Isometries;
Proof of Theorem \ref{main-2}}

Now we present a precise and detailed proof of Theorem
\ref{main-2}.
\bigskip

Denote by $\BB(p;r)=\{ z \in \CC^n \mid |z-p| < r \}$, the open
ball of radius $r$ centered at $p$ with respect to the Euclidean
distance on $\CC^n$.

Because of the $\mathcal{C}^2$ strong pseudoconvexity of
$\partial\Omega$ at the boundary orbit accumulation point $p$,
there exist a positive real number $\varepsilon$ and a
biholomorphic mapping $\Psi: U \to \BB (0; \varepsilon)$ such that
$\Psi(p)=0$,
\begin{multline*}
\Psi(\Omega \cap U) = %
 \{(z_1,\ldots,z_n) \in \BB (0; \varepsilon) \mid
\\
  \Re z_1 > |z_1|^2 + \ldots + |z_n|^2 %
   + E(z_1, \ldots, z_n)\ \} \, ,
\end{multline*}
and
\begin{multline*}
\Psi(\bo \cap U) = \{(z_1,\ldots,z_n) \in \BB (0; \varepsilon)
\mid
\\ \Re z_1 = |z_1|^2 + \ldots + |z_n|^2 + E( z_1, \ldots,
z_n)\ \},
\end{multline*}
where
$$
E(z_1, \ldots, z_n) = o(|z_1|^2 + \ldots + |z_n|^2).
$$
\medskip

Apply now the localization method by N. Sibony that uses only the
plurisubharmonic peak functions. (See \cite{SIB}, \cite{BER},
\cite{GAU}, \cite{BYGK} for instance, as well as \cite{ROY}.) It
follows from the hypothesis that every open neighborhood $U$ of
$p$ in $\CC^n$ admits an open set $V$ in $\CC^n$ such that $p \in
V \subset\subset U$ and
$$
\frac12 \ k_{\Omega \cap U} (z, \xi) \le k_{\Omega} (z, \xi)
$$
for every $z \in \Omega \cap V$ and every $\xi$ in the tangent
space $T_z \Omega~ (=\CC^n)$ of $\Omega$ at the point $z \in
\Omega$. It can also be arranged that
$$
\frac12 \ d_{\Omega \cap U} (x,y) \le d_\Omega (x,y)
$$
for every $x,y \in \Omega \cap V$. See \cite{BYGK} for instance
for this last inequality. This in particular implies
\medskip

\begin{quote} \noindent
\underbar{The localization property}: \sl For any open
neighborhood $W$ of $p$ and for any relatively compact subset $K$
of $\Omega$, there exists a positive integer $j_0$ such that
$\varphi_j(K) \subset \Omega \cap W$ whenever $j > j_0$. \rm
\end{quote}
\medskip

Notice that one can take $W$ such that $\Omega \cap W$ equipped
with the Kobayashi distance $d_{\Omega \cap W}$ is Cauchy
complete. Consequently the localization property, together with
the fact that $p$ is a boundary orbit accumulation point, implies
that $\Omega$ is complete Kobayashi hyperbolic.
\medskip

Next we apply Pinchuk's scaling method \cite{PIN}. Let
$$
\Psi \circ \varphi_j (q) \equiv (q_{1j}, \ldots, q_{nj})
$$
for $j=1,2,\ldots$. Fix $j$ for a moment. Choose $p_{1j} \in \CC$
such that
$$
(p_{1j}, q_{2j}, \ldots, q_{nj}) \in \Psi(\bo \cap U)
$$
and
$$
q_{1j} - p_{1j} > 0.
$$
Let $\zeta = A_j (z)$ for the complex affine map $A_j: \CC^n \to
\CC^n$ defined by
\begin{eqnarray*}
\zeta_1 & = & e^{i \theta_j}  (z_1 - p_{1j}) +
\sum_{k=2}^n c_{kj}(z_k - q_{kj}) \\
\zeta_2 & = & z_2 - q_{2j} \\
  & \vdots & \\
\zeta_n & = & z_n - q_{nj},
\end{eqnarray*}
where the real number $\theta_j$ and the complex numbers $c_{2j},
\ldots, c_{nj}$ are chosen so that the real hypersurface
$A_j\circ\Psi(\bo \cap U)$ is tangent to the real hyperplane
defined by the equation $\Re \zeta_1 = 0$.  It is important to
notice now, for the computation in the later part of this proof,
that
$$
\lim_{j\to\infty} e^{i \theta_j} = 1 ~\hbox{ and }~
\lim_{j\to\infty} c_{m j} = 0
$$
for every $m \in \{2,\ldots,n \}$.
\medskip

Then define $\Lambda_j : \CC^n \to \CC^n$ by
$$
\Lambda_j (z_1, \ldots, z_n) = \left (\frac{z_1}{\lambda_j},
\frac{z_2}{\sqrt{\lambda_j}}, \ldots, \frac{z_n}{\sqrt{\lambda_j}}
\right ) \, ,
$$
where $\lambda_j = q_{1j} - p_{1j}$.
\medskip

Exploit now the multi-variable Cayley transform
$$
\Phi (z_1, \ldots, z_n) = \left ( \frac{1-z_1}{1+z_1},
\frac{2z_2}{1+z_1}, \ldots, \frac{2z_n}{1+z_1} \right ) \, ,
$$
and consider the following sequences
$$
\sigma_j := \Phi \circ \Lambda_j \circ A_j \circ \Psi |_U \, ,
$$
and
$$
\tau_j := \sigma_j \circ \varphi_j \, ,
$$
for $j = 1,2,\ldots$.

Notice that each $\sigma_j$ maps $U$ into $\CC^n$.  It plays the
role of a holomorphic embedding of $\Omega \cap U$ into $\CC^n$.
On the other hand, the domain of definition of $\tau_j$ has to be
considered more carefully. Thanks to the {\it localization
property} above, for every compact subset $K$ of $\Omega$, there
exists a positive integer $j(K,U)$ such that $\tau_j$ maps $K$
into $\CC^n$ for every $j \ge j(K,U)$.  Thus, for each such $K$,
one is allowed to consider $\tau_j |_K : K \to \CC^n$ {\it only
for the indices $j$ with $j \ge j(K,U)$}.
\medskip

Now a direct calculation shows that, shrinking $U$ if necessary,
for every $\epsilon > 0$ there exists a positive integer $N$ such
that
$$
\BB (0; {1-\epsilon}) \subset \sigma_j (\Omega \cap W) \subset
\BB (0; {1+\epsilon})
$$
for every $j > N$.  Thus, replacing $N$ by $N + j(K,U)$, one may
conclude that
$$
\tau_j (K) \subset \BB (0;{1+\epsilon})
$$
for every $j > N$.
\medskip

Take now a sequence $\{K_\nu\}$ of relatively compact subsets of
$\Omega$ satisfying the following three conditions:
\begin{itemize}
\item[(\romannumeral1)] $K_\nu$ is a relatively compact, open
subset of $\Omega$ for each $\nu$;
\item[(\romannumeral2)] $\overline{K_\nu} \subset K_{\nu+1}$, for
$\nu = 1,2,\ldots$;
\item[(\romannumeral3)] $\displaystyle{\bigcup_{\nu=1}^\infty
K_\nu = \Omega}$.
\end{itemize}
Such a sequence $\{K_\nu\}$ is usually called a {\it (relatively)
compact exhaustion sequence} of $\Omega$.

Given a relatively compact exhaustion sequence $\{K_\nu\}$ of
$\Omega$, we consider the restricted sequences $\{\tau_{j,\nu} =
\tau_j |_{K_\nu} \mid j=1,2,\ldots \}$, for every
$\nu=1,2,\ldots$.
\medskip

Using these restricted sequences, we would like to establish:
\medskip

\begin{quote}
\noindent \bf Claim (\dag): \sl There exists a compact exhaustion
sequence $\{K_\nu\}$ such that the sequence $\tau_j$ admits a
subsequence that converges uniformly on every compact subset of
$\Omega$ to a Kobayashi isometry $\widehat\tau : \Omega \to \BB$
from the domain $\Omega$ onto the ball $\BB$. \rm
\end{quote}
\bigskip

Notice that the proof of Theorem \ref{main-2} is complete as soon
as this Claim is established.
\bigskip

Apply now Lemma \ref{lemma-1} to the domain $\Omega$. Let $\nu$ be
an arbitrarily chosen positive integer. Let $z_0 \in \Omega$. Let
the Kobayashi metric ball $B_{d_\Omega}(z_0; \mu\nu)$ play the
role of the subdomain $\Omega'$, where $\mu$ is an integer with
$\mu>5$. Then, for any $x, y \in B_{d_\Omega}(z_0,\nu)$, it holds
that
$$
d_{B_{d_\Omega}(z_;2\mu\nu)} (x,y) \le \frac1{\tanh (\mu\nu)} \
d_\Omega (x,y).
$$
\medskip

Exploiting the fact that $(\Omega, d_\Omega)$ is Cauchy complete,
we now choose the relatively compact exhaustion sequence
consisting of expanding Kobayashi metric balls:
$$
K_\nu \equiv B_{d_\Omega}(q;\nu).
$$
\medskip

Take $N>0$ such that $\varphi_j (q) \in V \cap \Omega$ and
$\varphi_j(K_{2\mu\nu}) \subset \Omega \cap U$ whenever $j > N$.
Enlarging $N$ if necessary, we may achieve also that $\sigma_j (\Omega \cap
U) \subset \BB(0;1+\epsilon)$ for every $j
> N$. Moreover, for any $x,y \in K_\nu$, one sees that:
\begin{eqnarray*}
d_{\BB (0; 1+\epsilon)} (\tau_j (x), \tau_j (y)) %
& \le & d_{\sigma_j(\Omega \cap U)}
       (\sigma_j \circ \varphi_j (x), \sigma_j \circ \varphi_j (y)) \\
& = & d_{\Omega \cap U} (\varphi_j (x), \varphi_j (y)) \\
& \le & d_{B_{d_\Omega} (\varphi_j(q),; 2\mu\nu))}
  (\varphi_j (x), \varphi_j (y)) \\
& \le & \frac{1}{\tanh (\mu\nu)}\ d_\Omega (\varphi_j(x), \varphi_j (y)) \\
& = & \frac{1}{\tanh (\mu\nu)}\ d_\Omega (x,y).
\end{eqnarray*}
As a summary, we have that
\begin{equation}
\label{estimate-1}
 d_{\BB (0, 1+\epsilon)} (\tau_j (x), \tau_j (y)) \le
 \frac1{\tanh (\mu\nu)} \ d_\Omega (x,y), \quad \forall x, y \in K_\nu.
\end{equation}

This estimate shows that the sequence $\{\tau_j\}$ is an
equi-continuous family on each $K_\nu$. Therefore one may extract
a subsequence that converges uniformly on every compact subset of
$\Omega$ to a continuous map $\widehat\tau : \Omega \to
\overline{\BB (0; 1+\epsilon)}$.
\medskip

Now the above estimate yields
\begin{equation}
\label{estimate-2} d_{\BB (0; 1+\epsilon)} (\widehat\tau (x),
\widehat\tau (y)) \le
 \frac1{\tanh (\mu\nu)} \ d_\Omega (x,y), \quad \forall x, y \in K_\nu.
\end{equation}
Since this estimate must hold for every $\epsilon>0$, one deduces
first that $\widehat\tau (K_\nu) \subset \overline{\BB(0;1)} =
\overline{\BB}$. But then, using the distance estimate above one
sees immediately that $\tau_j (K_\nu)$ for any $j$ is bounded away
from the boundary of $\BB$.  So $\widehat\tau (K_\nu) \subset \BB$
for every $\nu$. Consequently $\widehat\tau$ maps $\Omega$ into
$\BB$.
\medskip


Moreover,
$$
d_{\BB} (\widehat\tau (x), \widehat\tau (y)) \le \frac1{\tanh
(\mu\nu)} \ d_\Omega (x,y)), \qquad \forall x, y \in K_\nu.
$$
Letting $\mu \to \infty$, this last estimate turns into
\begin{equation}
\label{isom-1} d_{\BB} (\widehat\tau (x), \widehat\tau (y)) \le
d_\Omega (x,y)), \qquad \forall x, y \in K_\nu.
\end{equation}

Now let $x, y \in \Omega$ be fixed. Choose $0<\delta<1/2$ such
that $\widehat\tau (x), \widehat\tau(y) \in \BB(0; 1-2\delta)$.
Then there exists a positive integer $N_1$ such that
$$
\tau_j (x), \tau_j(y) \in \BB(0; 1-\delta)
$$
and
$$
\BB(0;1-\delta) \subset \sigma_j (\Omega \cap U)
$$
whenever $j>N_1$. Thus one obtains
\begin{eqnarray*}
d_\Omega (x,y) & = & d_\Omega (\varphi_j (x), \varphi_j (y)) \\
& \le & d_{\Omega \cap U} (\varphi_j (x), \varphi_j (y)) \\
& \le & d_{\sigma_j(\Omega \cap U)} (\sigma_j \circ \varphi_j (x),
  \sigma_j \circ \varphi_j (y)) \\
& = & d_{\sigma_j(\Omega \cap U)} (\tau_j (x), \tau_j (y)) \\
& \le & d_{\BB(0; 1-\delta)} (\tau_j (x), \tau_j (y)). \\
\end{eqnarray*}

Let $j$ tend to infinity first, and then let $\delta$ converge to
zero. Then one deduces that
\begin{equation}
\label{isom-2} d_\Omega (x,y) \le d_\BB (\widehat\tau (x),
\widehat\tau(y)).
\end{equation}

Combining (\ref{isom-1}) and (\ref{isom-2}), one sees that
$$
d_\Omega (x,y) = d_\BB (\widehat\tau (x), \widehat\tau(y)).
$$
Since $x$ and $y$ have been arbitrarily chosen points of $\Omega$,
it follows that $\widehat\tau : \Omega \to \BB$ preserves the
Kobayashi distance.
\bigskip

In order to complete the proof of the claim and hence Theorem
\ref{main-2}, it remains to show that $\widehat\tau : \Omega \to
\BB$ is surjective. Let $y \in \BB$. Then there exists $r$ with
$0<r<1$ such that $|y| < r$. Moreover, there exists $N_2>0$ such
that $\tau_j^{-1} (y) \in \Omega$ for every $j> N_2$.

Let $x_j = \tau_j^{-1} (y)$. Then it holds that $d_\Omega (q, x_j)
\le d_\BB (0; y) + 1$ for every $j > N_2$. Therefore a subsequence
$x_{j_k}$ of $x_j$ converges, say, to $\widehat x \in \Omega$. Now,
because of the uniform convergence of $\tau_j$ to $\widehat\tau$
on compacta, one immediately sees that
$$
\widehat\tau (\widehat x) = \widehat\tau (\lim_{j \to \infty}
(x_j)) = [\lim_{k \to \infty} \tau_k] (\lim_{j \to \infty} (x_j))
= \lim_{j\to\infty} \tau_j (x_j) = y.
$$
This shows that $\widehat\tau : \Omega \to \BB$ is surjective.
Consequently, the proof of Claim (\dag) follows.  The
proof of Theorem \ref{main-2} is now complete.
\endpf

\section{Complex Analyticity of the Kobayashi Isometry
$f:\Omega \to \BB^n$}

It now remains to establish Theorem \ref{main-main}. Incidentally,
it seems appropriate for us to pose the following naturally
arising question:

\begin{question} \sl
Let $n$ be a positive integer. Let $\Omega_1$ and $\Omega_2$ be
bounded domains in $\CC^n$ with $\mathcal{C}^2$ smooth, strongly
pseudoconvex boundaries, and let $f:\Omega_1 \to \Omega_2$ be a
homeomorphism that is an isometry with respect to the Kobayashi
distances.  Then, is $f$ or $\overline{f}$ necessarily
holomorphic?
\end{question}

We do not know the answer to this question at present; we show in
this paper that the answer is affirmative in case $\Omega_2 =
\BB^n$ and $\partial\Omega_1$ is $\mathcal{C}^{2,\epsilon}$
smooth.

\subsection{Burns-Krantz construction of Lempert discs
for strongly pseudoconvex domains}
Here we would like to explain how Burns and Krantz adapted
Lempert's analysis to the strongly pseudoconvex domains, as this
set of ideas is going to play an important role for our proof. In
what follows, let $\Omega \subseteq \CC^n$ be a bounded, strongly
pseudoconvex domain with ${\mathcal C}^{2,\epsilon}$ boundary.
\medskip

Let $p \in \partial\Omega$. Then by Burns and Krantz \cite{BUK},
there exist open neighborhoods $V$ and $U$ of $p$ such that $p \in
V \subset\subset U$ such that any $p' \in \partial\Omega \cap V$
admits a Lempert disc $\varphi:\overline\Delta \to
\overline{\Omega}$ such that $p' \in \varphi(\partial\Delta)$ and 
$\varphi(\overline{\Delta}) \subset U$. 
Perhaps the term Lempert disc needs to be clarified.  In fact
we will do more than that. We will quickly describe what Burns and
Krantz present in Proposition 4.3 and Lemma 4.4 of \cite{BUK}.
\medskip

Take the Fornaess embedding (\cite{FOR}) that embeds $\Omega$
holomorphically and properly into a strongly convex domain
$\Omega' \subset \CC^N$ with some $N>n$. This embedding map, say
$F$, is in fact smooth up to the boundary of $\Omega$ such that
$F:\overline{\Omega} \to \overline{\Omega'}$ is smooth, and also
$F(\partial\Omega) \subset \partial\Omega'$. Let $T_{F(p)}
(F(\Omega))$ be the tangent plane to $F(\Omega)$ at $F(p)$. This
is an $n$-dimensional complex affine space in $\CC^N$ and so we
may identify it with the standard $\CC^n$. Let $\Pi:\CC^N \to
T_{F(p)} (F(\Omega))$ be the orthogonal projection. Then $\Pi
\circ F:\Omega \to \CC^n$ is a holomorphic mapping, and
furthermore is a injective holomorphic mapping of $\Omega \cap U$
for some open neighborhood $U$ of $p$ in $\CC^n$.
\medskip

Denote by $F' := \Pi \circ F$. Since $F'(\Omega)$ is bounded,
there exists a sufficiently large ball $B$ such that $\Omega
\subset B$ and $p \in \partial\Omega \cap \partial B$. Slide $B$
slightly along the inward normal direction of $\partial F'(\Omega)
$ at $F'(p)$ and call it $B'$, so that the point $p$ is now
outside $B'$ and yet all the rest of $F'(\Omega)$ is within $B'$
except a small neighborhood $U'$ of $p$ satisfying $U' \subset U$.
Namely $F'(\Omega) \subset B' \cup U'$ with $U' \subset U$.
\medskip

Consider the convex hull the union of $F'(\Omega)$ and $B'$ and call it
$\Omega'$. This domain is convex and its boundary near $F(p)$
coincides the boundary of $F'(\Omega)$.  Note that the boundary of
$\Omega'$ is neither smooth nor strictly convex.  However it is
easy to modify $\Omega'$ slightly so that the newly modified domain
$\Omega''$ is strongly convex with $\mathcal{C}^6$ boundary, and
yet the boundary of $\Omega''$ coincides with the boundary of
$F'(\Omega)$ in a neighborhood of $F'(p)$.
\medskip

Then one considers Lempert discs, i.e., the holomorphic discs that
are isometric-and-geodesic embeddings of $\Delta$, for the domain
$\Omega''$. (See \cite{LEM2} for this.)  If one considers the
Lempert discs centered at a point sufficiently close to $p$ with
the direction at $p$ nearly complex parallel to a complex tangent
direction to the boundary of $\Omega''$ at $p$, then the image of
such discs will be within a small neighborhood of $p$. Moreover
such Lempert discs, say $h$, are also holomorphic geodesic
embeddings of the disc $\Delta$ into $\Omega$, in the sense that
$\check h := [F']^{-1} \circ h$ is the holomorphic geodesic
embedding of the unit disc $\Delta$ into $\Omega$. For further
details, see the above cited text in \cite{BUK}, especially
Proposition 4.3 and Lemma 4.4 therein.
\bigskip

For convenience, we shall call these discs {\it
Lempert-Burns-Krantz discs} for the strongly pseudoconvex domain
$\Omega$, or {\it LBK-discs} for short. Such discs exist at a
point sufficiently close to the boundary along the directions that
are approximately complex tangential to the boundary.

\subsection{Holomorphicity along the Lempert-Burns-Krantz discs}
Take now an LBK-disc $h:\Delta \to \Omega$
in the domain $\Omega$ such that $h^* d_\Omega =
d_\Delta$. Then we first present:

\begin{proposition}
For any continuous Kobayashi distance isometry $f: \Omega \to
\BB^n$, the composition $f\circ h:\Delta \to \BB^n$ is holomorphic
or conjugate-holomorphic.
\end{proposition}

\it Proof. \rm  Denote by $\widetilde h := f\circ h$. We give the
proof in two steps.
\medskip

{\bf Step 1: \it The mapping $\widetilde h$ is $C^\infty$ smooth.} We shall
first show that $\widetilde h$ is smooth at the origin. Take three
points $a, b$ and $c$ in the unit disc such that the Poincar\'e
geodesic triangle, say $T(a,b,c)$, with vertices at these three
points contains the origin in its interior. Fill $T(a,b,c)$ with
the geodesics from $a$ to points on the geodesic joining $b$ and
$c$.  Then obviously the origin is on one of these geodesics. Now,
let $m$ denote the foot of this geodesic. This procedure defines a
smooth diffeomorphism, say $h$, from a Euclidean triangle onto
$T(a,b,c)$, having two parameters: one is the time parameter of
each geodesic from the point $a$ to a point on the geodesic
joining $b$ and $c$, and the other is the parameter of the
geodesic joining $b$ and $c$.
\smallskip

Let $\widetilde a = \widetilde h (a), \widetilde b = \widetilde h
(b)$ and $\widetilde c = \widetilde h (c)$.  Let $\widetilde
T(\widetilde a,\widetilde b,\widetilde c)$ be the geodesic
triangle in $\BB^n$ with respect to the Poincar\'e metric, with
vertices at the three points $\widetilde a, \widetilde b$ and
$\widetilde c$. Again one may fill this triangle by Poincar\'e
geodesics of the ball, namely by the geodesics joining $\widetilde
a$ to the points on the geodesic joining $\widetilde b$ and
$\widetilde c$. This will again yield a smooth diffeomorphism from
a Euclidean triangle onto the filled triangle $\widetilde
T(\widetilde a,\widetilde b,\widetilde c)$. Since the Kobayashi
distance-balls are strongly convex for both $d_\Delta$ and
$d_{\BB^n}$, it follows that $\widetilde T(\widetilde a,\widetilde
b,\widetilde c) = \widetilde h (T(a,b,c))$. Moreover $\tilde h$
maps each geodesic to a corresponding geodesic with matching
speed. This shows that $\widetilde f$ is indeed smooth at the
origin. As the argument can be easily modified to prove the
smoothness of $\widetilde h$ at any point of the disc $\Delta$,
the map $\widetilde h:\Delta \to \BB$ is $\mathcal{C}^\infty$
smooth at every point.
\medskip

{\bf Acknowledgement:}  Notice that this argument can be used to give a
proof of the well-known theorem of Myers-Steenrod (\cite{MYS}).
This simple but elegant and powerful technique was conveyed to the
authors by Robert E. Greene in a private communication.  We
acknowledge with a great pleasure our indebtedness to him.
\bigskip

{\bf Step 2. \it The mapping $\widetilde h$ is holomorphic or conjugate
holomorphic.} Since $\widetilde h$ maps geodesics to geodesics, it
is a geodesic embedding.   Thus the surface $\widetilde
T(\widetilde a,\widetilde b,\widetilde c)$ has the maximal
holomorphic sectional curvature, and this can be realized only by
holomorphic sections in the ball. (Notice that the Kobayashi
metric coincides with the Poincar\'e metric in the unit ball, and
hence it is K\"ahler with negative constant holomorphic sectional
curvature.) Thus the tangent plane to the surface $\widetilde
T(\widetilde a,\widetilde b,\widetilde c)$ is complex. Since
$d\tilde h_* (T_* \Delta)$ is always a complex subspace in
$T_{\widetilde h (*)} \BB^n$, it follows by a standard argument
that $\widetilde h$ is either holomorphic or conjugate
holomorphic.
\endpf

\subsection{The Lempert map for strongly pseudoconvex domains}
Recall that our domain $\Omega''$ is a bounded strongly convex
domain with $\mathcal{C}^6$ boundary. Let $x \in \Omega''$. Then,
for an arbitrarily chosen $z \in \Omega''$ with $z \not= x$, there
exists a unique Lempert disc $h_{x,z}: \Delta \to $ such that
$h_{x,z} (0) = x$ and $h_{x,z} (\lambda) = z$ for some $\lambda$
with $0<\lambda<1$.  Then in \cite{LEM1, LEM2} Lempert defines
$\Phi(z) = \lambda h_{x,z}{}'(0)$ and shows that $\Phi:\Omega''
\to \CC^n$ extends to a $\mathcal{C}^2$ smooth map of the closure
$\overline{\Omega''}$. We shall call $x$ the {\it pivot} of the
map $\Phi$.

\noindent (Although Lempert was mainly interested in the
representation mapping $\widetilde\Phi (z) = \lambda
h_{x,z}{}'(0)/|h_{x,z}{}'(0)|$, which is nowadays known as the
Lempert representation map $\widetilde\Phi:\overline{\Omega''} \to
\overline{\BB^n}$, we remark here that both $\Phi$ and
$\widetilde\Phi$ are known to be $\mathcal{C}^2$ smooth up to the
boundary (\cite{LEM2}).) \bigskip

In the next subsection, we will see how to use this map $\Phi$ to show
that the Kobayashi isometry $f:\Omega \to \BB^n$ in question is
$\mathcal{C}^2$ up to the boundary.

\subsection{Smooth extension of Kobayashi isometry to the boundary}
Choose $p' \in S$ that is sufficiently close to $p$. Then let
$h:\Delta \to \Omega''$ be an LBK-disc with $h(0)=F'(p')$, as
mentioned above. Let us continue to use the notation $\check h :=
[F']^{-1} \circ h$. Then consider the M\"obius transformation
$\mu:\BB^n \to \BB^n$ which maps $f(p')$ to the origin.  The
preceding arguments imply that the composition $\widehat h:= \mu
\circ f \circ \check h:\Delta \to \BB^n$ defines a Lempert disc at
the origin for the unit ball $\BB^n$. Thus it is linear. Moreover,
it follows immediately that $\widehat h (\lambda) = \lambda
\widehat h'(0) = \lambda (\mu\circ f\circ \check h)'(0) =
d[\mu\circ f]_{p'} (\lambda \check h'(0))$.
\medskip

It is known that our $f$, a Kobayashi distance isometry, admits a
Lipschitz $1/2$ extension (\cite{HEN}). But with the strong
assumption that $\Omega$ is bounded strongly pseudoconvex with
$\mathcal{C}^6$ boundary, we shall prove the following:

\begin{proposition}
The Kobayashi isometry $f:\Omega \to \BB$ has a $\mathcal{C}^2$
extension to the boundary.  More precisely, there exists an open
neighborhood $W$ of $\partial\Omega$ in $\CC^n$ such that
$f:\overline{\Omega}\cap W \to \overline{\BB}$ is $\mathcal{C}^2$
smooth.
\end{proposition}

\it Proof. \rm
Notice first that the Kobayashi isometry $f$ as well as $f^{-1}$
are locally Lipschitz with respect to the Euclidean metric, as the
Kobayashi distance generates the same topology as the Euclidean
distance (\cite{BAR}). Therefore the set
$$
S := \{x \in \Omega \mid df_x \hbox{ and } d[f^{-1}]_{f(x)} \hbox{
exist } \}
$$
is a subset of full measure, i.e., the Lebesgue measure of $S$ is
the same as the Lebesgue measure of $\Omega$.  In particular, $S$
is dense in $\Omega$.
\medskip

Now let $\Upsilon (h(\lambda)) := [d[\mu \circ f]_{p'}]^{-1} \circ
\mu \circ f (\check h (\lambda))$ for $\lambda \in \Delta$. Since
$\mu \circ f \circ \check h (\lambda) = d[\mu\circ f]_{p'}
(\lambda \check h'(0))$, it follows that $\Upsilon$ coincides, at
every point on the image $h(\Delta)$, with the aforementioned map
$\Phi:\Omega'' \to \CC^n$ with its pivot at $p'$.  To be precise
for any $\zeta \in h(\Delta)$, let $h(\lambda)=\zeta$. Then
$\Upsilon (\zeta) = \lambda h'(0)$.
\medskip

Altogether, one sees that the mapping $[d[\mu \circ f]_p]^{-1} \circ \mu
\circ f \circ [F']^{-1}$ coincides with $\Upsilon$ in a small
conical neighborhood (with apex at $p$) of $h(\partial\Delta)$
filled by the LBK-discs at $p'$; it follows that $f$ is
$\mathcal{C}^2$ smooth in an open neighborhood of $\check h
(\partial\Delta)$. Now it is easy to observe that this gives rise
to the $\mathcal{C}^2$ smoothness of $f$ in an open neighborhood
of $\partial\Omega$ as desired.
\endpf
\bigskip

\subsection{The Kobayashi isometry is CR on the boundary}
We now present

\begin{proposition} \label{CR}
The restriction of the extension of $f$ to $\partial\Omega$ into
$\CC^n$ is a CR function (or an anti-CR function).
\end{proposition}

\it Proof. \rm  Let $p \in \partial\Omega$ and let $\overline{L}
\in T_p^{0,1} \partial\Omega$.  Regard this vector field as a
derivation operator on the Euclidean space $\CC^n$. Then, for every
$\epsilon > 0$, there exists an $r \in (0,\epsilon)$ and $q \in
\Omega$ with $|p-q|=r$ such that we may find a Lempert disc
$\varphi:\Delta \to \Omega$ satisfying
$$
\overline{L}_q f = \frac{\partial}{\partial \bar\zeta}\Big|_0
f\circ \varphi.
$$
\medskip

Since $f \circ \varphi$ is holomorphic (replace it by $\bar f
\circ \varphi$ if necessary), one immediately sees that
$\overline{L}_q f = 0$  Since $f$ is $\mathcal{C}^2$ up to the
boundary, letting $r$ tend to zero, one obtains the assertion.
\endpf

\subsection{Analyticity of the Kobayashi isometry -- proof of
Theorem \ref{main-main}}
Finally we are ready to present:
\bigskip

\it Proof of Theorem \ref{main-main}. \rm Start with Proposition
\ref{CR} which we just proved. Recall that $f$ restricted to
$\partial\Omega$ is a $C^2$ smooth CR map from $\partial\Omega$ to
$\partial\BB$.

It is also a local diffeomorphism.  Notice that $f(\partial\Omega)$
is compact and relatively open. Therefore $f(\partial\Omega) = \partial\BB$.
This implies that $f$ is a covering map.  However $\partial\Omega$
is simply connected if $n>1$, being topologically a sphere.  Hence
$f:\partial\Omega \to \BB$ is a $C^2$ diffeomorphism.

Now apply the Bochner-Hartogs theorem.  The mapping $f$ extends to a
holomorphic mapping, say $\widehat f$ of $\Omega$ into $\BB$. Now
restrict $f$ to a Lempert-Burns-Krantz disc.  This has the
same value as $\widehat f$ at any point of the disc. Therefore,
$f$ and $\widehat f$ must coincide at every point of the LBK-disc.
Altogether, the map $f$ itself is holomorphic in $W \cap \Omega$ for some
open neighborhood $W$ of $\partial\Omega$.  

Then one may ask whether $f=\widehat f$ on $\Omega$. They do coincide indeed. This can be seen as follows. Since $f:\Omega \to \BB$ is a Kobayashi distance isometry, $\Omega$ itself has the property that any two points in it must have one and the only one shortest distance realizing curve joining them. Take any such curve $\gamma:[0,\ell] \to \Omega$ with $\gamma(0), \gamma(\ell) \in W \cap \Omega$. Then one observes the following: 
\begin{itemize}
\item $d_\Omega (\gamma(s), \gamma(t)) = s-t, \hbox{ whenever } 0\le t \le s \le \ell$.
\item $f \circ \gamma$ is the unique distance realizing curve joining $f(\gamma(0))$ and $f(\gamma(\ell))$. 
\item $f(\gamma(0))= \widehat f(\gamma(0))$ and $f(\gamma(\ell))=\widehat f(\gamma(\ell))$.
\end{itemize} 
Now for every $t \in [0,\ell]$ notice that
$$
d_\Omega (f\circ\gamma(0), \widehat f\circ \gamma (t)) = d_\Omega (\widehat f\circ\gamma(0), \widehat f\circ \gamma (t)) \le d_\Omega (\gamma(0),\gamma(t)) = t,
$$
and likewise
$$
d_\Omega (\widehat f\circ \gamma (t), f\circ\gamma(\ell)) \le \ell-t.
$$
Now by triangle inequality and the uniqueness of the (shortest) distance realizing curve joining two points, one sees immidiately that the inequalities above are equalities and that 
$$
\widehat f (\gamma(t)) = f (\gamma(t))
$$
for every $t \in [0,\ell]$. It is now immediately deduced that $f=\widehat f$ on $\Omega$. In particular $f$ (and hence $\widehat f$) is a bijective holomorphic mapping of $\Omega$ onto the ball $\BB$. This finishes the
proof. \endpf
\bigskip

\it Note: \rm The authors would like to thank H.\ Seshadri for asking whether $f$ can be shown directly to coincide with $\widehat f$; we clarified it changing the end of the proof slightly. 

\section{Concluding Remarks}

The original theorem of Bun Wong [WON], and variants
by Rosay [ROS] and others, has proved to embody a
powerful and far-reaching set of ideas.  In particular,
it was consideration of this insight that led Greene
and Krantz ([GRK2], [GRK3], [GRK4]) to formulate
the principle that the Levi geometry of a boundary
orbit accumulation point will determine the global
geometry of the domain.  This in turn has led
to the Greene-Krantz conjecture:  that a boundary
orbit accumulation point for a smoothly bounded domain
must in fact be of finite type in the sense of
Kohn-Catlin-D'Angelo.

The result of these studies has been a profound and
fruitful development in geometric analysis.  We wish
that the present contribution will lead to further insights.

\vspace{1in}

Kang-Tae Kim

Department of Mathematics

Pohang University of Science and Technology

Pohang 790-784 The Republic of Korea (South)

\it e-mail: \tt kimkt@postech.ac.kr \rm
\bigskip   \\

Steven G. Krantz


American Institute of Mathematics

360 Portage Avenue

Palo Alto, CA 94306-2244  U.S.A.

\it e-mail: \tt skrantz@aimath.org \rm

\end{document}